\documentclass[a4paper]{jpconf}
\usepackage[utf8]{inputenc}
\usepackage[english]{babel}

\usepackage{calc}
\usepackage{amssymb}
\usepackage{amstext}
\usepackage{amsmath}
\usepackage{amsthm}
\usepackage{mathrsfs}
\usepackage{dsfont}
\usepackage{graphicx}
\graphicspath{{./images/}{./}}
\usepackage{pslatex}
\usepackage{url}
\usepackage[ruled,vlined]{algorithm2e}
\usepackage{xspace}
\usepackage[caption=false]{subfig}
\usepackage{placeins}

\newtheorem{thm}{Theorem}

\newcommand{\norm}[1]{\left\Vert #1\right\Vert}
\newcommand{\abs}[1]{\left\vert #1\right\vert}
\newcommand{\dotproduct}[2]{\ensuremath{\langle #1, #2 \rangle}\xspace}
\newcommand{\p}{\partial}

\newcommand{\R}{ \mathbb{R} }
\newcommand{\Z}{ \mathbb{Z} }
\newcommand{\N}{ \mathbb{N} }

\newcommand{\LL}{ \mathbb{L} }

\newcommand{\argmin}{\mathop{\mathrm{arg\,min}}}

\usepackage{color}
\begin{document}
\title{Image restoration using sparse approximations of spatially varying blur operators in the wavelet domain}

\author{Paul Escande$^1$, Pierre Weiss$^1$, and Fran\c{c}ois Malgouyres$^2$}
\address{$^1$ ITAV-USR3505, Universit\'e de Toulouse, CNRS, Toulouse, France.}
\address{$^2$ IMT-UMR5219, Universit\'e de Toulouse, CNRS, Toulouse, France.}

\ead{paul.escande@gmail.com, pierre.armand.weiss@gmail.com, francois.malgouyres@math.univ-toulouse.fr}

\begin{abstract}
Restoration of images degraded by spatially varying blurs is an issue of increasing importance in the context of photography, satellite or microscopy imaging.
One of the main difficulty to solve this problem comes from the huge dimensions of the blur matrix.
It prevents the use of naive approaches for performing matrix-vector multiplications.
In this paper, we propose to approximate the blur operator by a matrix sparse in the wavelet domain. 
We justify this approach from a mathematical point of view and investigate the approximation quality numerically.
We finish by showing that the sparsity pattern of the matrix can be pre-defined, which is central in tasks such as blind deconvolution. 
\end{abstract}

\section{Introduction}

The problem of image restoration in the presence of spatially varying blurs drew a considerable attention in the context of satellite imaging with the Hubble space telescope \cite{nagy1998restoring}. It is now becoming increasingly important in various imaging modalities such as 3D fluorescence microscopy \cite{RR_SBH_11}.
In this paper, we aim at providing new computational tools to tackle this problem. We consider a blurring operator $H$ in $\Omega \subset \R^d$, defined for any $u \in \LL^2(\Omega)$ as the following integral operator,
\begin{equation}
	\label{eq:inte_operator}
	\forall x \in \Omega, \quad Hu(x) = \int_{y \in \Omega} K(x,y) u(y) dy,
\end{equation}

where $\Omega \subset \R^d$ is the image domain. The function $K$ is a spatially varying kernel defining the Point Spread Function (PSF) at each location $x \in \Omega$. 

The naive computation of a product $Hu$ is simply obtained by direct discretization of \eqref{eq:inte_operator} and costs $\mathcal{O}\left(N^{2d}\right)$ arithmetic operations for an image of dimension $d$ with $N$ pixels in each dimension. It is unsuitable for large scale problems encountered in imaging. Previous works attempted to reduce this computational cost by approximating $H$ using either tensor products or piecewise convolutions \cite{zhang2007gaussian,nagy1998restoring, hansen2006deblurring}. Unfortunately, these strategies do not have a property that is highly desirable in numerical analysis: the ability to approximate the original operator with an arbitrary precision. In practice, we observed that they are often too coarse for complex situations encountered in microscopy imaging.



In a recent work \cite{escande2012spatially}, we proposed to approximate $H$ by an operator diagonal in the wavelet domain. 
However, this approximation was shown to be too coarse for some practical applications. Our main contribution in this paper consists of showing that the diagonal approximation can be replaced by a \textit{sparse} approximation and lead to arbitrarily close approximations. 
More precisely, we show that $\|H - \Psi S \Psi^*\|$, can be made arbitrarily small if $\Psi^*$ denotes a forward wavelet transform and $S$ is a sparse matrix that contains only $O(N^d)$ non zero coefficients instead of $N^{2d}$. As a practical consequence, matrix-vector products using the proposed approach costs $\mathcal{O}\left(N^{d}\right)$ arithmetic operations, which is doable even for very large images. We also show that the sparsity pattern of $S$ can be known a priori, which is a highly desirable property in the case of blind deconvolution, where the operator should be inferred from the blurred images.

In section \ref{sec:theoretical}, we recall the theoretical results motivating this approximation. 
Section \ref{sec:threshold}, contains preliminary restoration results in order to illustrate theoretical results and to identify a number of coefficients offering a compromise between the storage costs and the restoration quality. In this section, the kept coefficients are obtained by simple thresholding and can have arbitrary coordinates. 
In section \ref{sec:patterns}, we finally show that the sparsity pattern of the matrix $\Theta$ can be pre-defined and that this approach leads to almost equivalent results to a thresholding.

\section{Theoretical motivations} \label{sec:theoretical}

	Sparse approximations of integral operators have been theoretically analyzed in \cite{BCR, coifman1997wavelets}. Surprisingly this approach was never applied to the approximation of spatially varying blur operators. The closest application found in the literature deals with foveation \cite{CMY-Foveation}, but it is not adapted to a large class of kernels.
	
	We define an orthonormal wavelet basis of $\LL^2(\R)$ as the family of functions:
	
	\[
		\left\{ \psi_{j,m} \right\}_{j,m \in \Z}.
	\]
	Each $\psi_{j,m}$ is a dilated and translated of the mother wavelet $\psi$,
	\[
		\psi_{j,m}(t) = \sqrt{2^{-j}} \psi\left( 2^{-j} t -m \right).
	\]

	Let us recall a typical result that motivates the proposed approach. We stick to the one-dimensional case for the ease of exposition. 
Since $H$ is a linear operator in a Hilbert space, it can be written as $H = \Psi \Theta \Psi^*$, where $\Theta: l^2 \rightarrow l^2$ is the matrix representation of the blur operator in the wavelet domain. Going from an infinite dimensional setting to a finite dimensional one can be performed by using projectors on linear subspaces of the multi-resolution analysis. Matrix $\Theta$ is characterized by the coefficients:
	\[
		\theta_{j,m,k,n} := \left( \dotproduct{H \psi_{j,m}}{\psi_{k,n}} \right), \qquad \forall j,m,k,n \in \Z.
	\]
	
	
	\begin{thm}[Decay of $\theta_{j,m,k,n}$ -- \cite{BCR}]\label{thm:beylkin}
	Supposing that:	
		\begin{itemize}
			\item the wavelets are compactly supported with $\mathop{supp}(\psi_{j,m})=I_{j,m}$,
			\item the wavelets have $M$ vanishing moments,
			\item the operator $H$ belongs to the class of Calderon-Zygmund operators, meaning that $K$ satisfies:
				\begin{equation}
					\exists C_M > 0, \text{ such that }, \forall x, y \in \Omega, \quad \left\{
					\begin{split}
						\abs{K(x,y)} & \leq \frac{1}{\abs{x-y}}, \\
						\abs{\p_x^M K(x,y)} + \abs{\p_y^M K(x,y)} & \leq \frac{C_M}{\abs{x-y}^{1+M}},
					\end{split} 
					\right.
				\end{equation}
		\end{itemize}
		\bigskip
		then the coefficients $\theta_{j,m,k,n}$ satisfy the following estimate:
		\begin{equation}
			\abs{\theta_{j,m,k,n} } \leq C_M 2^{-\frac{\abs{j-k}}{2}} \left( \frac{ 2^{\min(j,k)} }{ \textrm{dist}(I_{j,m}, I_{k,n})} \right)^{M+1},
		\end{equation}
		with $C_M$ a constant depending on $M$ and the wavelets, $\textrm{dist}(I_{j,m}, I_{k,n})$ denotes the distance between the two supports of the wavelets $\psi_{j,m}$ and $\psi_{k,n}$.
		
		Moreover, if the kernel $K$ is compactly supported, then for sufficiently large $\abs{m-n}$,
		\begin{equation}
			\theta_{j,m,k,n} = 0.
		\end{equation}
	\end{thm}

\begin{figure}[h!]
\begin{center}
\begin{minipage}{0.45\textwidth}
	\centering
		\includegraphics[height=0.8\textwidth]{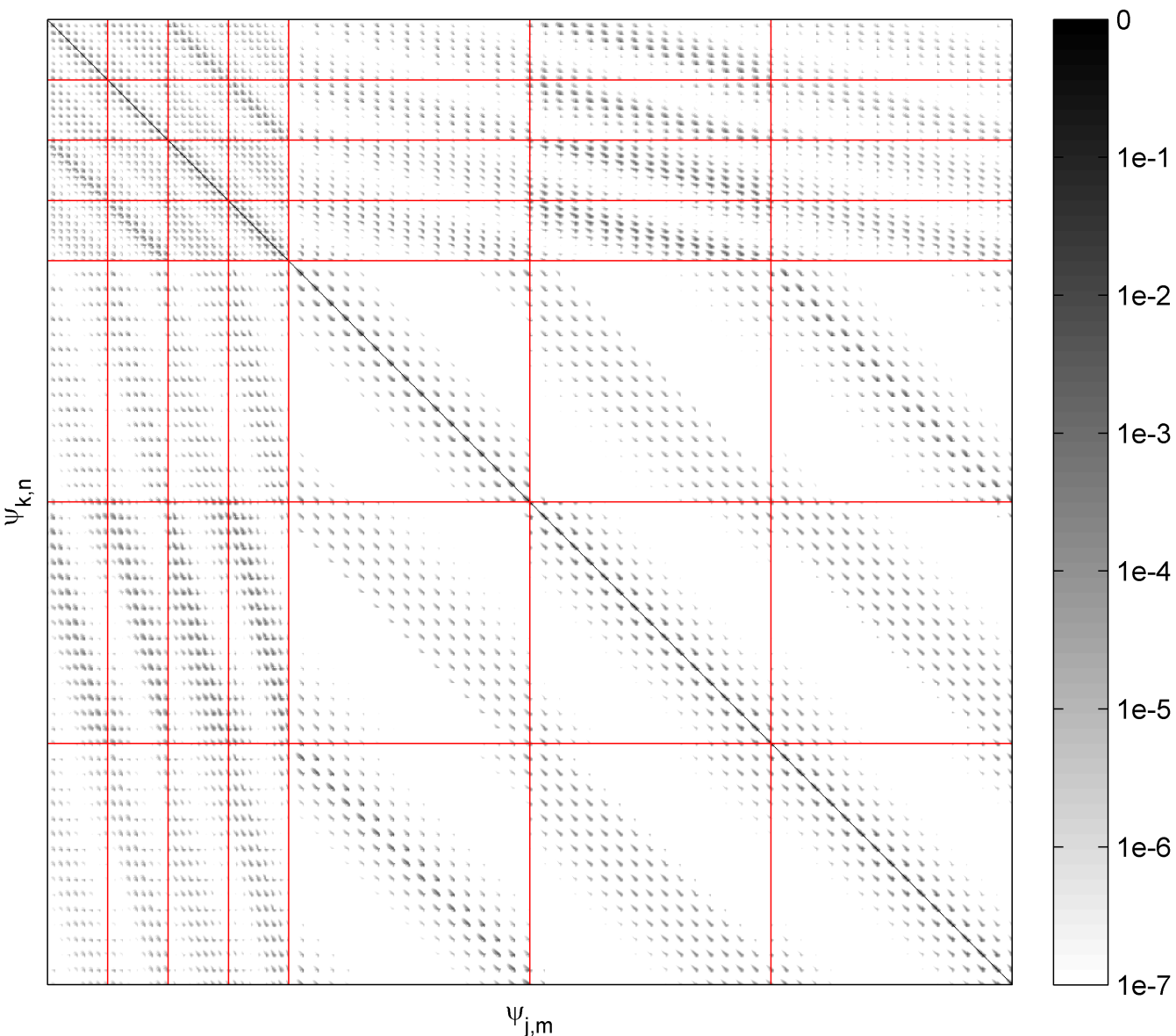}
		\caption{An example of a $\Theta$ matrix in a $\log_{10}$ scale for a $256\times 256$ images. The matrix contains about $4.3$ billion coefficients.} \label{fig:logTheta}	
\end{minipage}\hspace{0.02\textwidth}%
\begin{minipage}{0.45\textwidth}
\centering
\includegraphics[height=0.65\textwidth]{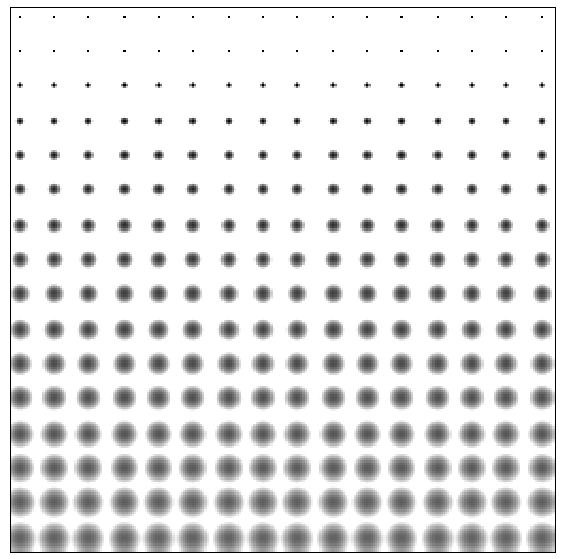}
\caption{The point spread function of the blur used to generate the matrix on the left. The kernel is a Gaussian with equal variances increasing in the vertical direction.}
\end{minipage} 
\end{center}
\end{figure}

	This results ensures that the coefficients away from the diagonal fastly decay to zero. A more precise analysis shows that $\Theta$ contains only $\mathcal{O}\left(N\epsilon^{-1/(M+1)}\right)$ above $\epsilon$ \cite{BCR, meyer1992wavelets}. 
	Figure \ref{fig:logTheta} shows the matrix $\Theta$ associated to a typical spatially varying blur operator in log scale. 
	It is readily seen that the decrease of the coefficients away from the diagonal is extremely fast.

	In the next section, we reconstruct degraded images from the knowledge of various thresholded $\Theta$ matrices. 

	\section{Restoration results}
	
	\subsection{Restoration using thresholded matrices} \label{sec:threshold}
	
\begin{figure}[htp]
	\begin{center}
	\subfloat[Original Image \label{fig:u}]{\includegraphics[width=0.22\textwidth]{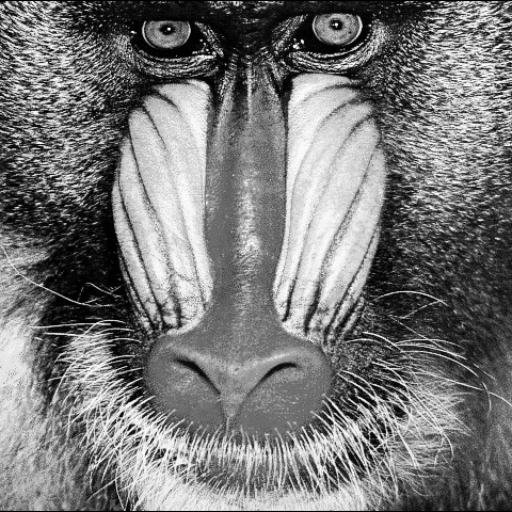}} $\quad$
	\subfloat[Blurred Image]{\includegraphics[width=0.22\textwidth]{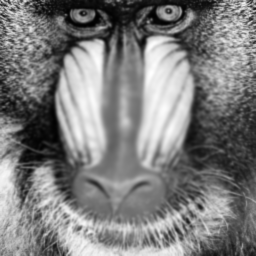}} $\quad$
	\subfloat[Degraded Image $\sigma = 2.10^{-2}$, $SNR = 16.54$dB.\label{fig:v}]{\includegraphics[width=0.22\textwidth]{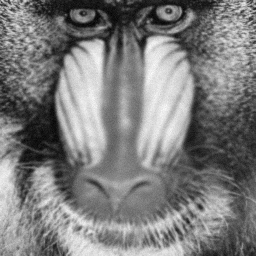}}$\quad$
	\subfloat[PSF\label{fig:PSF}]{\includegraphics[width=0.22\textwidth]{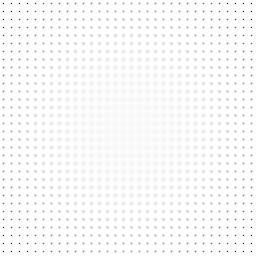}}
		\end{center}
	\caption{Images involved in the deblurring problem}
\end{figure}

	In a first experiment, the $\Theta$ matrix is thresholded in order to keep a number of coefficients being a multiple of $N^{d}$ (the number of pixels of the image). That is, we define the thresholded matrix $\Theta_{T_k}$ by zeroing all but the $T_k = k N^d$ largest coefficients of $\Theta$. This experiment allows the identification of an admissible number of coefficients offering a proper compromise between the storage costs and the image restoration quality. To achieve this experiment, a full $\Theta$ matrix is numerically computed and then successively thresholded for various $k \in \N^*$. Thereafter, the image is restored assuming the following classical degradation model:
	\[
		v = Hu + \eta, \quad \eta \sim \mathcal{N}(0, \sigma^2 I_{n^d}),
	\]
	where $v$ is the degraded image observed (see Figure \ref{fig:v}), $u$ is the image to restore (see Figure \ref{fig:u}), and $H$ is a blurring operator corresponding to the PSF displayed in Figure \ref{fig:PSF}. A standard TV-L2 optimization problem is solved to restore the image $u$:
	
\begin{equation}
\textrm{Find \ } u^*\in \argmin_{u \in  \R^{n^d}, \norm{\Psi \Theta_{T_k} \Psi^* u - v}_2^2 \leq \sigma^2 n} TV(u),
\label{eq:prob_deblurring}
\end{equation}
where $TV(u)$ is the isotropic total variation of $u$.
	
\begin{figure}[htp]
	\begin{center}
	\subfloat[Non-zeros coefficient for $k=1$]{\includegraphics[width=0.29\textwidth]{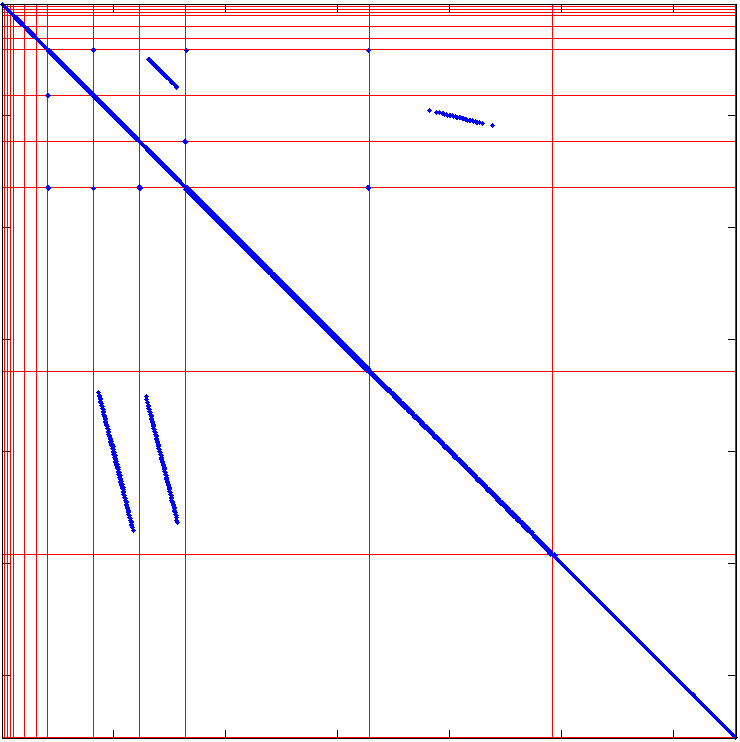}} $\quad$
	\subfloat[Non-zeros coefficient for $k=20$]{\includegraphics[width=0.29\textwidth]{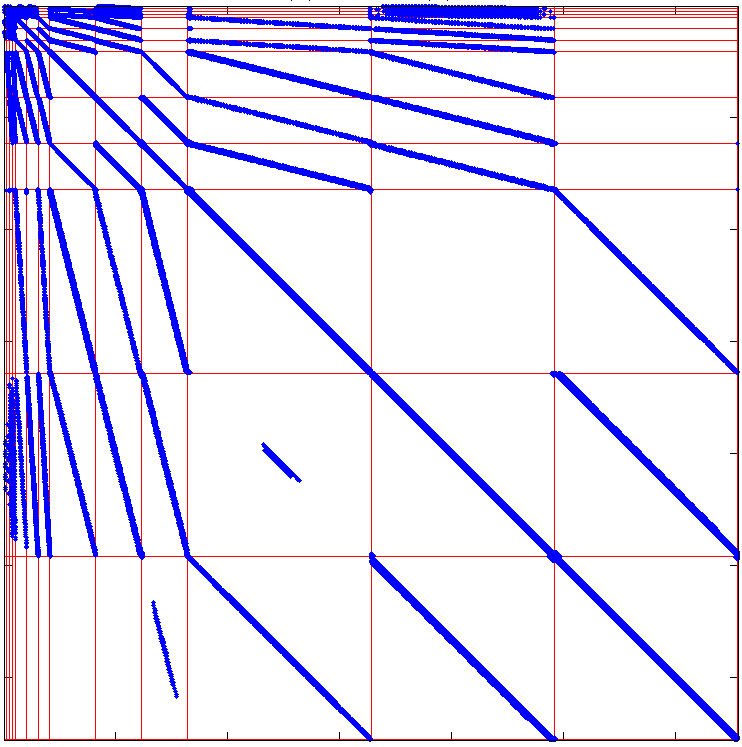}} $\quad$
	\subfloat[Full $\Theta$ matrix]{\includegraphics[width=0.29\textwidth]{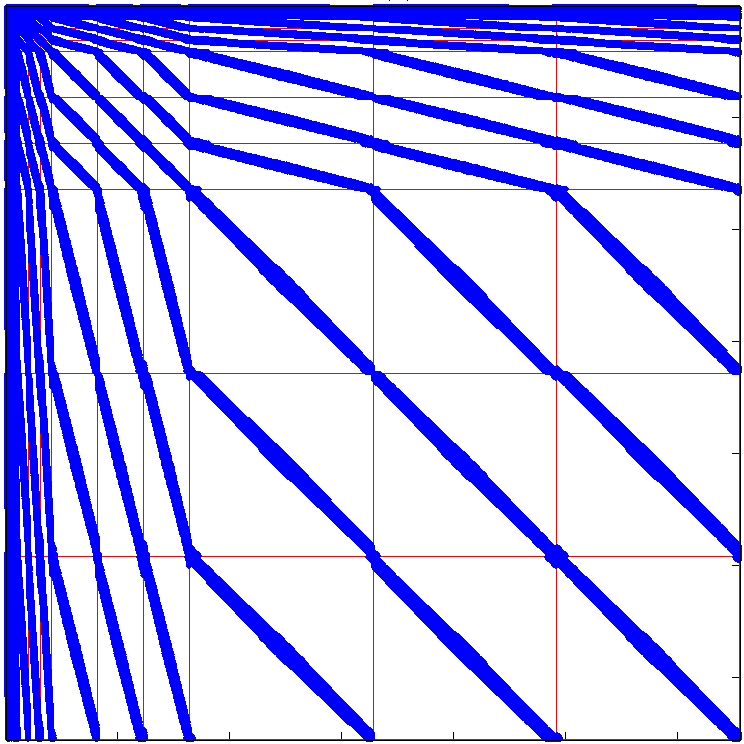}}
	
	\subfloat[Deblurred image for $k=1$. $SNR = 17.12$dB, Speed-up = 22. \label{fig:thresh_1}]{\includegraphics[width=0.29\textwidth]{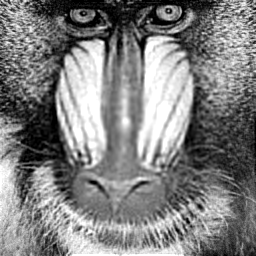}}$\quad$
	\subfloat[Deblurred image for $k=20$. $SNR = 18.69$dB, Speed-up = 17. \label{fig:thresh_20}]{\includegraphics[width=0.29\textwidth]{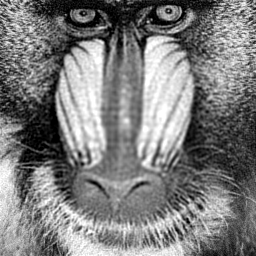}}$\quad$
	\subfloat[Deblurred image with the full $\Theta$, $k \simeq 1550$, $SNR = 18.93$dB. \label{fig:thresh_full}]{\includegraphics[width=0.29\textwidth]{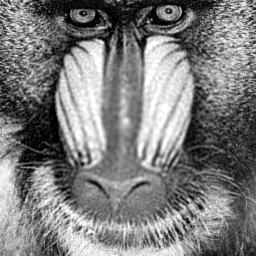}}
		\end{center}
	\caption{Deblurred images using thresholded matrices} \label{fig:thresh}
\end{figure}

	Figure \ref{fig:thresh} displays restored images for various values of $k$ and the  sparsity pattern of $\Theta_{T_k}$. As expected, considering more coefficients improves the deblurring results. However, in the presence of noise, the two images restored using the full matrix and the thresholded matrix with $k=20$, provide very similar results. Therefore, it is possible to store only a small amount of coefficients and obtain results close to the ones obtained having a full (but unrealistic) knowledge of the blur matrix. Furthermore, these results are obtained drastically faster, see the speed-up Figure \ref{fig:thresh}.
	
	\paragraph{Remark:} it is near impossible to solve the deblurring problem with an exact operator due to the numerical burden of computing matrix-vector products in large dimensions. We thus chose a spatially varying blur that is sparse in the wavelet domain in order to be able to make comparisons between an exact computation and an approximated one. Speed up times might thus be much higher for other kinds of blurs.
	
%

	\subsection{Restoration using pre-defined sparsity patterns} \label{sec:patterns}

	We are currently developping strategies to solve blind deconvolution problems using the proposed approximation. 
In this setting, it is not possible to threshold the exact (unknown) matrix to obtain a sparsity pattern. 
The knowledge of such a structure is primordial to approximate $K$ using prior information. 
Our aim in this section is to show that sparsity patterns can be pre-defined, allowing convincing restoration results. 

Our main observation is that for a given wavelet $\psi_{j,m}$, only its neighbours in the same scale, the finer scale and the coarser one lead to significant correlation coefficients $ \dotproduct{H \psi_{j,m}}{\psi_{k,n}}$. This concept of neighbours is also supported by Theorem \ref{thm:beylkin}. Thus, sparsity patterns can be defined using notions of intra and inter-scale neighbourhoods. 

\begin{figure}[h!]
\begin{center}
\begin{minipage}{0.45\textwidth}
\centering
\includegraphics[width=0.7\textwidth]{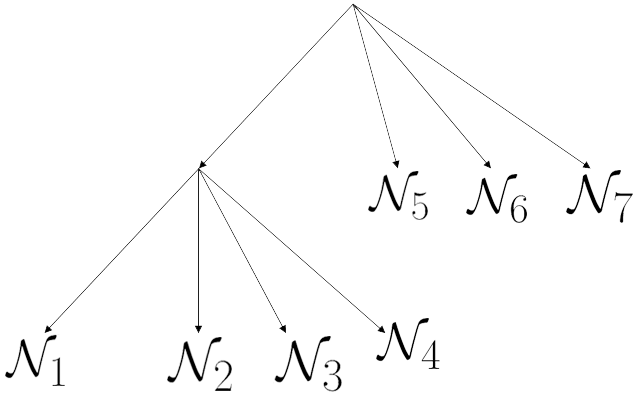}
\caption{Quadtree corresponding to a wavelet transform of depth 2} \label{fig:tree}
\end{minipage}\hspace{0.02\textwidth}%
\begin{minipage}{0.45\textwidth}
\centering
\includegraphics[width=0.6\textwidth]{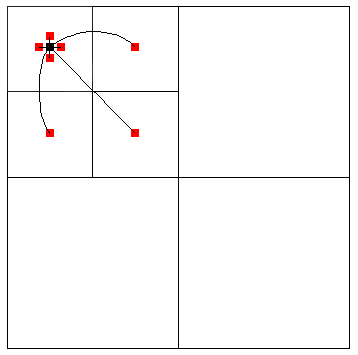}
\caption{Illustration of the neighbourhood $\mathcal{N}_1$. In black, the given wavelet at the center of the image, and in red all the preserved neighbouring wavelets.} \label{fig:N1}
\end{minipage} 
\end{center}
\end{figure}

For each node of the wavelet decomposition quadtree, a neighbourhood is defined. It describes which correlation coefficients shall be preserved to generate the sparsity pattern. Figure \ref{fig:tree} illustrates this idea. In this example, we consider a wavelet transform of depth 2. A neighbourhood $\mathcal{N}_i$ is associated to each sub-band. For instance $\mathcal{N}_1$ can be represented as follows:
\begin{equation*}
	\mathcal{N}_1 = \begin{pmatrix}
		2 & 2 & 2 & 2 & 2 & 2 & 2 & 2 \\
		 l & l & l & l & l & h & v & d \\
		 0 & -1 & 1 & 0 & 0 & 0 & 0 & 0 \\
		 0 & 0 & 0 & -1 & 1 & 0 & 0 & 0 \\
	\end{pmatrix}
	\left|\begin{array}{l}
	  \textrm{scale} \\
	  \textrm{oriention} \\
	  \textrm{vertical translation} \\
	  \textrm{horizontal translation}
	\end{array} \right.
\end{equation*}
where $l$ stands for the low frequency wavelet, $h, v, d$ for the horizontal, vertical and diagonal orientations respectively.
Figure \ref{fig:N1} illustrates this neighbourhood $\mathcal{N}_1$ for a given wavelet at the center of the image.

\begin{figure}[htp]
	\begin{center}
	\subfloat[Non-zero coefficients of pattern 1]{\includegraphics[width=0.29\textwidth]{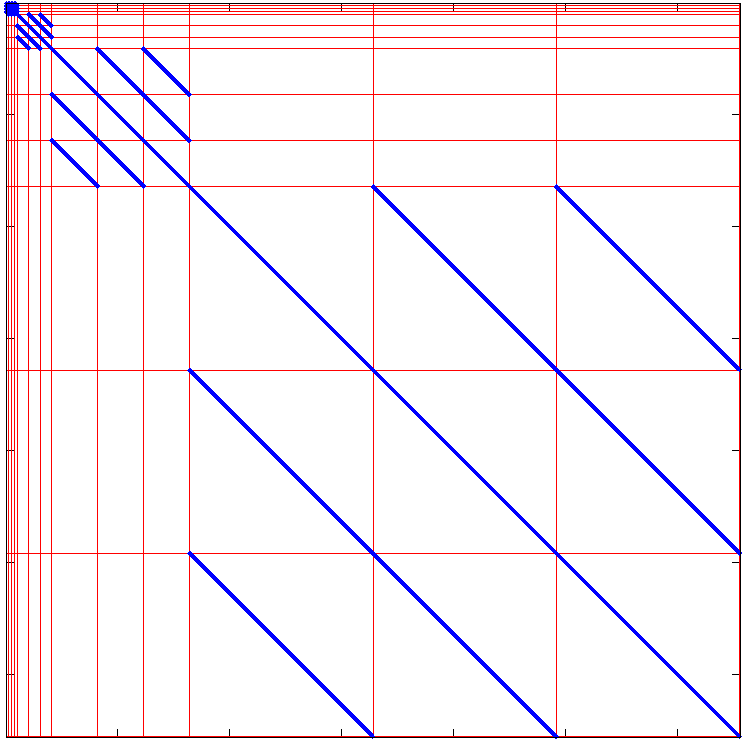}} $\quad$
	\subfloat[Non-zero coefficients of pattern 2]{\includegraphics[width=0.29\textwidth]{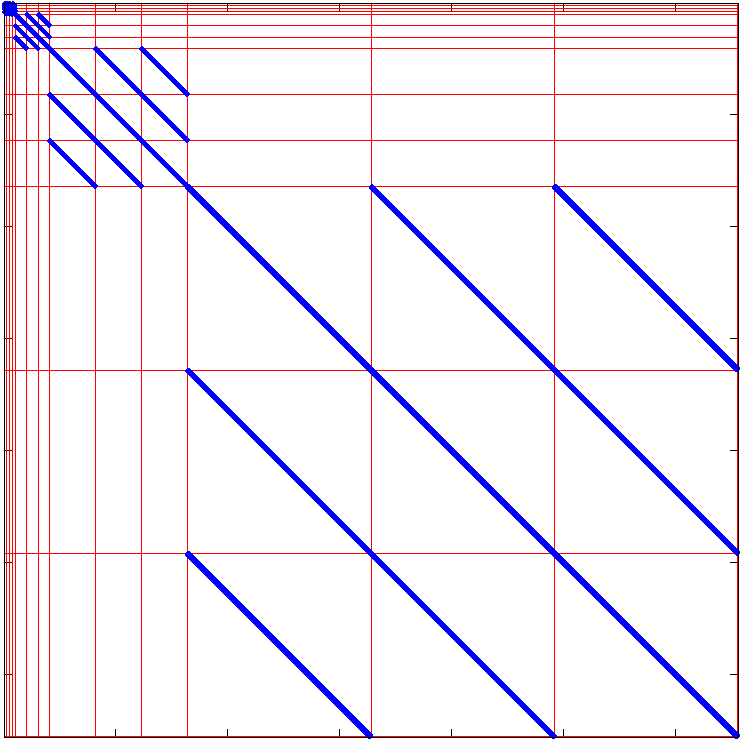}} 
	
	\subfloat[Deblurred image with pattern 1, $ k \simeq 3$, $SNR = 16.74$dB.]{\includegraphics[width=0.29\textwidth]{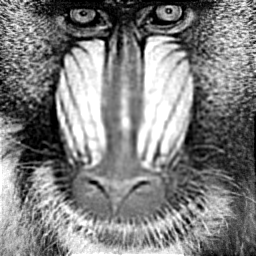}}$\quad$
	\subfloat[Deblurred image with pattern 2, $k \simeq 15$, $SNR = 18.30$dB.]{\includegraphics[width=0.29\textwidth]{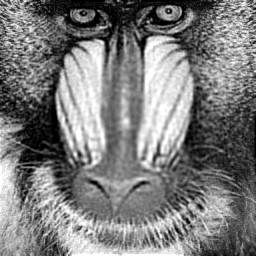}}
		\end{center}
	\caption{Deblurred images using matrices defined by sparsity patterns} \label{fig:pattern}
\end{figure}

In order to validate this notion of neighbourhood, we considered two different scenarii. The first one consists of setting for all $i$, $\displaystyle \mathcal{N}_i = \begin{pmatrix}
		\text{same scale} \\
		\text{all orientations} \\
		0 \\
		0
	\end{pmatrix}.$
The second one consists of setting for all $i$,

\[
	\mathcal{N}_i = \begin{pmatrix}
		\text{same scale} & \text{same scale} & \text{same scale} & \text{same scale} & \text{same scale} \\
		\text{all orientations} & \text{all orientations} & \text{all orientations} & \text{all orientations} & \text{all orientations} \\
		0 & -1 & 1 & 0 & 0\\
		0 & 0 & 0 & -1 & 1 \\
	\end{pmatrix},
\]
Figure \ref{fig:pattern} displays restored images using sparsity patterns generated from these two neighbourhoods.
This experiment shows that sparsity patterns defined from the correlations between neighbouring wavelets are relevant for our aim: approximating and interpolating the blurring operator.

\section*{Conclusion}

In this paper we investigated the use of sparse matrices in the wavelet domain in order to approximate spatially varying blur operators. 
We observed that very large compression factors still allow to obtain near optimal reconstruction results.
Compared to previously proposed approach, this technique allows to approximate the exact operator with an arbitrary precision.
It is thus suited to noiseless and noisy problems.
We also showed that pre-established sparsity patterns can be used to efficiently approximate the blurring integral operator.
This property is central in order to tackle blind deconvolution problems. We are currently working on this topic.

\ack
This project was funded by ANR SPHIM3D.

\section*{References}
\nocite{*}
\bibliography{Biblio_Escande_ICPRAM_2012}

\end{document}